\documentclass[11pt]{article}
\usepackage[cp1251]{inputenc}
\usepackage[ukrainian, russian, english]{babel}

\usepackage{amsfonts, amssymb,mathrsfs}                                              %

\begin{document}
\selectlanguage{russian} \thispagestyle{empty}
 \pagestyle{myheadings}



\setcounter{page}{1}

\thispagestyle{empty}

\begin{center}

{\bf  ON SOME INEQUALITIES OF CHEBYSHEV TYPE }

\vskip 5mm

A{\small NDRIY} {L.}  {S}{\small HIDLICH}, S{\small TANISLAV} {O.}  {C}{\small HAICHENKO}

\vskip 5mm

{\small {\it Abstract.} We obtain some new inequalities  of Chebyshev Type.}

\vskip 5mm

{\bf 1. Introduction.}

\end{center}

\vskip 5mm

 Let $f$, $g$: $[a,b]\to {\mathbb R}$ be integrable functions, both increasing or both decreasing.
Further, let $p\,{:}\,[a,b]\to {\mathbb R}_0^+$ be an integrable
function. Then (see, for example, [1, Chap.~IX])
    \begin{equation}\label{p.1}
\int_a^b p(x) f(x)g(x)dx\ge
  \int_a^bp(x)f(x)dx\int_a^b p(x) g(x)dx\bigg(\int_a^bp(x)dx\bigg)^{-1}.
 \end{equation}
If one of the functions  $f$ or $g$ is nonincreasing and the other
nondecreasing the reversed inequality is true, i.e.,
    \begin{equation}\label{p.2}
\int_a^b p(x) f(x)g(x)dx\le
  \int_a^bp(x)f(x)dx\int_a^b p(x) g(x)dx\bigg(\int_a^bp(x)dx\bigg)^{-1}.
 \end{equation}
Inequalities  (\ref{p.1}) and (\ref{p.2}) are known as
Chebyshev's inequalities. These inequalities were obtained by P.L.~Chebyshev [2, 3] and they attracted great interest of the researchers.  So, a lot of analogues and generalizations
of inequalities  (\ref{p.1}) and (\ref{p.2}) is  known.  In particular, these
results can be found in Chapter IX of the book
[1] by D.S.~Mitrinovi\'c,
J.E.~Pe\v cari\'c and A.\,M.~Fink which trace
completely the historical and chronological developments of
Chebyshev's and related inequalities (see also [4, 5]). Also we would like to recommend the article of H.P.~Heinig and L.~Maligranda [6], where one can found a lot of important results on Chebyshev's inequalities for strongly increasing functions, positive convex and concave functions as well as on Chebyshev's inequalities in Banach function spaces and symmetric spaces.

In [7], these investigations were developed in the following direction: the author found necessary and sufficient conditions on the function $g$: $[a,b]\to {\mathbb R}_0^+$ and $p$: $[a,b]\to {\mathbb R}^+$ such that for any monotone function $f$: $[a,b]\to {\mathbb R}_0^+$ the relations
     \begin{equation}\label{p.5}
\int\limits_a^b p(x) f(x)g(x)dx\ge
 \bigg( \int\limits_a^bp^r(x)f^r(x)dx\bigg)^{1/r}\int\limits_a^b p(x) g(x)dx
  \bigg(\int\limits_a^bp^r(x)dx\bigg)^{-1/r}
 \end{equation}
and
      \begin{equation}\label{p.6}
\int\limits_a^b p(x) f(x)g(x)dx\le
 \bigg( \int\limits_a^bp^r(x)f^r(x)dx\bigg)^{1/r}\int\limits_a^b p(x) g(x)dx
  \bigg(\int\limits_a^bp^r(x)dx\bigg)^{-1/r}
 \end{equation}
hold with $r$ being an arbitrary positive number.

In this paper we continue the study of the inequalities of the type (\ref{p.1})--(\ref{p.6}), namely,  we obtain the following assertions:

T{\small HEOREM } 1. {\it Assume that  $g$: $[a,b]\to {\mathbb R}_0^+$ and $p$: $[a,b]\to {\mathbb R}^+$ are integrable functions such that the product $p\cdot g$ is also integrable on $[a,b]$ function. Let $f\,{:}\,[a,b]\to {\mathbb R}_0^+$ be a nonicreasing function.  Then for any convex function $M$: $[0,\infty)\to {\mathbb R}$ such that $M(0)=0$,  the following inequality is true:
\begin{equation}\label{p.5}
\int_a^b p(x) g(x) M(f(x)) dx \le  \sup\limits_{s\in (a,b]} \left\{M\bigg(\frac{\int_a^b p(x)f(x)dx }{\int_a^s p(x)dx}\bigg)\int_a^s p(x) g(x) dx \right\},
\end{equation}
and for any concave function $M$: $[0,\infty)\to {\mathbb R}$ such that $M(0)=0$,  the following inequality is true:
\begin{equation}\label{p.51}
\int_a^b p(x) g(x) M(f(x)) dx \ge  \inf\limits_{s\in (a,b]} \left\{M\bigg(\frac{\int_a^b p(x)f(x)dx }{\int_a^s p(x)dx}\bigg)\int_a^s p(x) g(x) dx \right\}.
\end{equation}
Furthermore, if the function $f(x)\equiv c$, $c\ge 0$, then relations $(\ref{p.5})$ and $(\ref{p.51})$ are equalities. }

Putting  $M(t)=t^{1/r}$, $r>0$, from Theorem 1 we obtain the following corollaries:

C{\small OROLLARY} 1.  {\it Let $r\in (0,1]$, and  let $g$:
$[a,b]\to {\mathbb R}_0^+$ and $p$: $[a,b]\to {\mathbb R}^+$  be
integrable functions such that for all $s\in (a,b]$,
 \begin{equation}\label{p.701}
\frac{\int_a^s p(x) g(x) dx }{\bigg(\int_a^s p(x)dx\bigg)^{1/r}} \le \frac{\int_a^b p(x) g(x) dx }{\bigg(\int_a^b p(x)dx\bigg)^{1/r}}.
\end{equation}
 Then for any nonicreasing function $f$: $[a,b]\to {\mathbb R}_0^+$,
\begin{equation}\label{p.511}
\int_a^b p(x) g(x) f(x) dx \le \bigg(\int_a^b p(x)f^r(x)dx\bigg)^{1/r}
\frac{\int_a^b p(x) g(x) dx }{\bigg(\int_a^b p(x)dx\bigg)^{1/r}}.
\end{equation}}

C{\small OROLLARY} 2.  {\it Let $r\ge 1$, and  let $g$:
$[a,b]\to {\mathbb R}_0^+$ and $p$: $[a,b]\to {\mathbb R}^+$  be
integrable functions such that for all $s\in (a,b]$,
 \begin{equation}\label{p.7011}
\frac{\int_a^s p(x) g(x) dx }{\bigg(\int_a^s p(x)dx\bigg)^{1/r}} \ge \frac{\int_a^b p(x) g(x) dx }{\bigg(\int_a^b p(x)dx\bigg)^{1/r}}.
\end{equation}
 Then for any nonicreasing function $f$: $[a,b]\to {\mathbb R}_0^+$,
\begin{equation}\label{p.5111}
\int_a^b p(x) g(x) f(x) dx \ge \bigg(\int_a^b p(x)f^r(x)dx\bigg)^{1/r}
\frac{\int_a^b p(x) g(x) dx }{\bigg(\int_a^b p(x)dx\bigg)^{1/r}}.
\end{equation}}

\vskip -5mm

If in corollaries 1 and 2, we put $r=1$, then we see that relations  (\ref{p.511}) and (\ref{p.5111})
 are the   Chebyshev's classical  inequalities (\ref{p.1}) and (\ref{p.2}). Furthermore, it should be noted that conditions on the functions $p$ and $g$ of the form  (\ref{p.701}) and (\ref{p.7011}) for validity of inequalities (\ref{p.1}) and (\ref{p.2}) were considered in the papers [7] and [8].

In the case, where the function $M(f(x))$ is nonincreasing and the function $g$ is nondecreasing (or nonincreasing), we can apply the  Chebyshev's classical inequalities to the integral $\int_a^b p(x) g(x) M(f(x)) dx $ on the left-hand side of relations (\ref{p.5}) (or (\ref{p.51})). Respectively, we obtain
\begin{equation}\label{p.52}
\int_a^b p(x) g(x) M(f(x)) dx \le \frac{\int_a^b p(x)M(f(x))dx }{\int_a^b p(x)dx}\int_a^b p(x) g(x) dx
\end{equation}
and
\begin{equation}\label{p.53}
\int_a^b p(x) g(x) M(f(x)) dx \ge \frac{\int_a^b p(x)M(f(x))dx }{\int_a^b p(x)dx}\int_a^b p(x) g(x) dx.
\end{equation}
Furthermore, if exact upper (or lower) bound on the right--hand side of  (\ref{p.5}) (or (\ref{p.51})) is realized for $s=b$, then from relations  (\ref{p.5}) and (\ref{p.51}) we get \vskip -5mm
\begin{equation}\label{p.54}
\int_a^b p(x) g(x) M(f(x)) dx \le  M\bigg(\frac{\int_a^b p(x)f(x)dx }{\int_a^b p(x)dx}\bigg)\int_a^b p(x) g(x) dx,
\end{equation}
and
\begin{equation}\label{p.55}
\int_a^b p(x) g(x) M(f(x)) dx \ge M\bigg(\frac{\int_a^b p(x)f(x)dx }{\int_a^b p(x)dx}\bigg)\int_a^b p(x) g(x) dx.
\end{equation}
Here, it should be note that by virtue of Jensen's inequality (see, for example, [1, Chap.~I]), estimations (\ref{p.54}) and (\ref{p.55}) of the integral $\int_a^b p(x) g(x) M(f(x)) dx $ are more precisely, than  estimations (\ref{p.52}) and (\ref{p.53}).

R{\small EMARK.} {In the case, where the function $f$ does not decrease, inequalities (\ref{p.5}) and (\ref{p.51}) have the similar form, but in these  inequalities, all the integrals of the kind  $\int\limits_a^s (\cdot)$  should be replaced by the integrals of the kind $\int\limits_s^b (\cdot)$. }

\vskip 2mm

\centerline{\bf 2. Discrete analogue of Theorem 1.}

 L{\small EMMA} 1.  {\it Assume that $a=\{a_k\}_{k=1}^m$,  $b=\{b_k\}_{k=1}^m$ and $p=\{p_k\}_{k=1}^m$, $m\in {\mathbb N}$ are nonnegative number sequences such that $a_1\ge a_2\ge\ldots a_m$ and $p_k>0$. Then for any convex function $M$: $[0,\infty)\to {\mathbb R}$ such that $M(0)=0$,  the following inequality is true:
 \begin{equation}\label{p.7}
 \sum_{k=1}^m p_k b_k M(a_k )\le \max\limits_{s\in [1,m]}\left\{   M\bigg(\frac {\sum_{k=1}^m p_k a_k }{\sum_{k=1}^s p_k}\bigg) \sum_{k=1}^s p_k b_k\right\},
 \end{equation}
 and for any concave function $M$: $[0,\infty)\to {\mathbb R}$ such that $M(0)=0$,  the following inequality is true:
 \begin{equation}\label{206}
 \sum_{k=1}^m p_k b_k M(a_k )\ge \min\limits_{s\in [1,m]}\left\{   M\bigg(\frac {\sum_{k=1}^m p_k a_k }{\sum_{k=1}^s p_k}\bigg) \sum_{k=1}^s p_k b_k\right\}.
 \end{equation}
 Furthermore, if the sequence $a_k\equiv c$, $c\ge 0$, then relations $(\ref{p.7})$ and $(\ref{206})$ are equalities. }

{\it Proof.} Consider the case, where the function  $M$ is convex (in the case, where the function  $M$ is concave, the proof is similar). Let us prove by the induction on $m$ the proposition that for any convex function $M$: $[0,\infty)\to {\mathbb R}$ such that $M(0)=0$, inequality (\ref{p.7}) holds.

The case $m=1$ is obvious.

Also consider the case $m=2$.

Put
 \begin{equation}\label{p.8}
c=p_1 a_1+p_2 a_2,\quad x_0=p_1 a_1,\quad \alpha_k=p_kb_k,\quad \beta_k=p_k^{-1},\ \  k=1,2,
 \end{equation}
and consider on the interval $[0,c]$ the function
 \begin{equation}\label{p.9}
h(x)=\alpha_1 M(\beta_1  x)+\alpha_2 M(\beta_2(c-x)).
 \end{equation}
 Due to convexity of the function  $M(t)$, the function $h(x)$  is also convex on $[0,c].$ Hence, this function attains its maximum  value on any
interval $[x_1,x_2] \subseteq [0,c]$ at one of its endpoints. Thus
  \begin{equation}\label{105}
 h(x)\le \max\{h(x_1),h(x_2)\}\quad \forall x\in [x_1,x_2].
 \end{equation}
Setting $x_1:=\beta_2 c (\beta_1+\beta_2)^{-1}$ and $x_2:=c$,  we see that the number $x_0$  (by virtue of monotonicity of the sequence $a$) belongs to the interval $[x_1,x_2]$.

Therefore, in view of  relations (\ref{p.8})--(\ref{105}) and of the equality $M(0)=0$, we get
$$
 \sum_{k=1}^2 p_k  b_kM(a_k)=h(x_0)\le \max\{h(x_1),h(x_2)\}
 =
 $$
 $$
 =\max\bigg\{M\Big(\frac{p_1a_2+p_2a_2}{p_1+p_2}\Big)(p_1b_2+p_2b_2),M\Big(\frac{p_1a_2+p_2a_2}{p_1}\Big)p_1b_1\bigg\}.
$$
Hence, for $m=2$, inequality (\ref{p.7}) holds.

Now, assume that for $m=n-1\ge 1$, the proposition is true.

Let us show that for  $m=n$,  it is also true. Let us use notations (\ref{p.8}) and  consider on the interval $[0,c]$   the function $h(x)$ of the form as in (\ref{p.9}).  Setting $x_1:=\beta_2 c (\beta_1+\beta_2)^{-1}$ and $x_2:=c-a_3/\beta_2$,  we see that the number $x_0$ (by virtue of monotonicity of the sequence $a$) belongs to the interval $[x_1,x_2]$. Thus  in view of relations  (\ref{p.8})--(\ref{105}),
\begin{equation}\label{200}
 \sum_{k=1}^n p_k  b_kM(a_k)=h(x_0)+ \sum_{k=3}^n p_k  b_kM(a_k)\le \max\{h(x_1),h(x_2)\}+ \sum_{k=3}^n p_k  b_kM(a_k).
 \end{equation}

Further, in the case, where $h(x_1)\ge h(x_2)$, we set
\begin{equation}\label{201}
 p{\,}'_k=\left\{\matrix{p_1+p_2,\ \hfill k=1,\cr
 \ \ p_{k+1},\ \ \  \hfill k=\overline{2,m-1};}\right.\quad
  b'_k=\left\{\matrix{(p_1 b_1+p_2b_2)/(p_1+p_2 ),\ \hfill k=1,\cr
 \quad\quad b_{k+1},\ \ \ \hfill k=\overline{2,m-1};}\right.
 \quad \ \ \
 \end{equation}
\begin{equation}\label{202}
a'_k=\left\{\matrix{(p_1 a_1+p_2 a_2)/(p_1+p_2),\ \hfill k=1,\cr
 \ \ \ \quad\quad a_{k+1},\ \ \  \hfill k=\overline{2,m-1}.}\right.
\end{equation}
Then by virtue of  (\ref{200}), we conclude that the following relation is true:
\begin{equation}\label{203}
\sum_{k=1}^m p_k  b_kM(a_k)
\le \sum_{k=1}^{m-1} p'_k  b'_k M(a'_k).
\end{equation}
In the case, where $h(x_1)< h(x_2)$, relation (\ref{203}) holds for the sequences $a'$, $b'$ and $p'$ of the form:
\begin{equation}\label{204}
 p{\,}'_k=\left\{\matrix{ \ \ p_{1},\ \ \  \hfill k=1,\cr  \ \ p_2+p_3,\ \hfill k=2,\cr
 \ \ p_{k+1},\ \ \  \hfill k=\overline{3,m-1};}\right.\quad
  b'_k=\left\{\matrix{\ \ b_{1},\ \ \  \hfill k=1,\cr \ \ (p_2 b_2+p_3b_3)(p_2+p_3)^{-1},\ \hfill k=2,\cr
 \ \  b_{k+1},\ \ \ \hfill k=\overline{3,m-1};}\right.
 \quad \ \ \
\end{equation}
 \begin{equation}\label{205}
a'_k=\left\{\matrix{\ \ (p_1a_1+p_2a_2-p_2a_3)/p_1,\ \ \  \hfill k=1,\cr
 \ \  a_{k+1}.\ \ \  \hfill k=\overline{2,m-1},}\right.
  \end{equation}
The sum on the right-hand side of (\ref{203}) contains $n-1$ items. Furthermore,  in both cases, for the sequences $a'$, $b'$ and $p'$, the induction assumption is satisfied. Thus, taking into account (\ref{201})--(\ref{205}), we obtain the necessary estimate (\ref{p.7}):
$$
\sum_{k=1}^n p_k  b_kM(a_k)
\le\sum_{k=1}^{n-1} p'_k  b'_k M(a'_k)\le \sup\limits_{s\in [1,n-1]}\left\{   M\left(\frac {\sum_{k=1}^{n-1} p'_k a'_k }{\sum_{k=1}^s p'_k}\right) \sum_{k=1}^s p'_k b'_k\right\}\le
$$
$$
\le \sup\limits_{s\in [1,n]}\left\{   M\left(\frac {\sum_{k=1}^{n} p_k a_k }{\sum_{k=1}^s p_k}\right) \sum_{k=1}^s p_k b_k\right\}.
$$
Lemma is proved.

\vskip 2mm

\centerline{\bf 3. Proof of Theorem 1.}

 Consider the case, where the function  $M$ is convex (in the case, where the function  $M$ is concave, the proof is similar). First, let us verify  that inequality (\ref{p.5}) holds  for any function $f$ such that
 for a certain $m\in {\mathbb N}$,
 $$
  f(x)=a_k,\quad x\in [l_{k-1},l_k),\ k=1,2,\ldots, m,
 $$
where $a_1>a_2>\ldots>a_{m}\ge 0$ and $a=l_0<l_1<\ldots<l_{m}=b$.

For any  $k=1,2,\ldots,m$, we put
 $$
 p_k=\int\limits_{l_{k-1}}^{l_k}p(x)dx,\quad
 b_k=\int\limits_{l_{k-1}}^{l_k}p(x)g(x)dx
 \bigg(\int\limits_{l_{k-1}}^{l_k}p(x)dx\bigg)^{-1}.
 $$
Then by virtue of Lemma 1, we get  (\ref{p.5}):
$$
\int_a^b p(x) g(x) M(f(x)) dx = \sum\limits_{k=1}^{m} \int\limits_{l_{k-1}}^{l_k} p(x) g(x) M(f(x))dx=
 \sum_{k=1}^m p_k b_k M(a_k )\le
 $$
 $$
 \le \sup\limits_{s\in [1,m]\cap {\mathbb N}}\left\{   M\bigg(\frac {\sum_{k=1}^m p_k a_k }{\sum_{k=1}^s p_k}\bigg) \sum_{k=1}^s p_k b_k\right\}=\sup\limits_{s\in [1,m]\cap {\mathbb N}}\left\{   M\bigg(\frac {\int_a^b p(x)f(x)dx }{\int_a^{l_s} p(x)dx}\bigg) \int_a^{l_s} p(x) g(x)dx\right\}\le
 $$
 $$
\le \sup\limits_{s\in (a,b]}\left\{   M\bigg(\frac {\int_a^b p(x)f(x)dx }{\int_a^{s} p(x)dx}\bigg) \int_a^{s} p(x) g(x)~dx\right\}.
$$

Let us prove the validity of inequality (\ref{p.5}) in general case. First, note that if the functions $M$ and $f$ satisfy the conditions of Theorem 1, then there exists the number $n_0=n_{0}(M,f)\in {\mathbb N}$ such that for any $n>n_0$ and for all $x\in [a;b]$, the inequality  $|M(f(x))|< n$ holds.

For any $n>n_0$, consider the system of points $l_0^{(n)}< l_1^{(n)}<\ldots< l_m^{(n)}=b$, defined in the following way:  we put $l_0^{(n)}:=a$  and for any $k\in [1;m]\cap \mathbb{N}$ the value $l_k^{(n)}$ is a greatest positive number such that $l_k^{(n)}>l_{k-1}^{(n)}$ and for all $x\in [l_{k-1}^{(n)};k_{k}^{(n)})$ the following relation is true:
$$
    |M(f(l_{k-1}^{(n)}))-M(f(x))|\le \frac 1n.
$$
By virtue of the conditions on the function $M$ and $f$, this system of points always exist and  $m\le 2n^2$.

Further, consider the functions  $f_n=f_n(x)$ such that
\begin{equation}\label{l.32q}
    f_n(x)\equiv \lim_{t \to l_{k}^{(n)}-} f(t),\quad \mbox{\rm for all}\
    x\in [l_{k-1}^{(n)}; l_{k}^{(n)}),\ k=1,2,\ldots, m.
\end{equation}
We see that the inequality  $|M(f(x))-M(f_n(x))|\le \frac 1n$ holds  for   all $n>n_0$ and $x\in [a,b]$. Due to integrability on $[a,b]$ of the product   $p(x)g(x)$, the values
$$
  \int_{a}^b p(x)g(x) [M(f(x))-M(f_n(x))]~dx
$$
converge to zero as $n\to\infty$.  Furthermore, for any  $n>n_0$, the function $f_n(x)$ is nonincreasing  and it takes finitely
many values on $[a,b]$. Hence, this function satisfies the conditions of the proposition proved above.

Thus, in view of (\ref{l.32q}) and continuity of  the function $M$, we conclude that for any  $\varepsilon>0$ and for all sufficiently great  $n$ ($n>n_1(\varepsilon)$),
 $$
  \int\limits_{a}^b p(x) g(x)M(f(x)) dx=  \int\limits_{a}^b p(x)g(x) M(f_n(x)) dx+
    \int\limits_{a}^b p(x)g(x) (M(f(x))-M(f_n(x))) dx\le
$$
$$
     \le \sup\limits_{s\in (a;b]}\left\{   M\bigg(\frac {\int_a^b p(x)f_n(x)dx }{\int_a^{s} p(x)dx}\bigg) \int_a^{s} p(x) g(x)dx\right\}+\frac{\varepsilon}{2}
    \le\sup\limits_{s\in (a;b]}\left\{   M\bigg(\frac {\int_a^b p(x)f(x)dx }{\int_a^{s} p(x)dx}\bigg) \int_a^{s} p(x) g(x)dx\right\}+\varepsilon.
 $$
Hence, relation  (\ref{p.5}) is true.

Analyzing the proof of Theorem 1, we see that the similar statement is also true in the case, where $b=\infty$.

T{\small HEOREM } $1'$.  {\it Assume that  $g$: $[a,b]\to {\mathbb R}_0^+$ and $p$: $[a,b]\to {\mathbb R}^+$  (where $b\in (a,\infty]$) are integrable functions such that the product $p\cdot g$ is also integrable on $[a,b]$ function. Let also $f\,{:}\,[a,b]\to {\mathbb R}_0^+$ be a nonicreasing function.  Then for any convex (or concave) function $M$: $[0,\infty)\to {\mathbb R}$ such that $M(0)=0$,  inequality  $(\ref{p.5})$  (or inequality  $(\ref{p.51})$) is true. }

Analogically, one can obtain the statement, similar to Lemma 1, in the case, where $n=\infty$.

 L{\small EMMA}  $1'$.  {\it Let  $a=\{a_k\}_{k=1}^\infty$,  $b=\{b_k\}_{k=1}^\infty$ and $p=\{p_k\}_{k=1}^\infty$ be nonnegative number sequences such that $a_1\ge a_2\ge\ldots$, $p_k>0$ and the series $\sum_{k=1}^\infty p_k b_k$ is convergent. Then for any convex function $M$: $[0,\infty)\to {\mathbb R}$ such that $M(0)=0$,  the following inequality is true:
$$
 \sum_{k=1}^\infty p_k b_k M(a_k )\le \sup\limits_{s\in [1,\infty)}\left\{   M\bigg(\frac {\sum_{k=1}^\infty p_k a_k }{\sum_{k=1}^s p_k}\bigg) \sum_{k=1}^s p_k b_k\right\},
 \eqno(\ref{p.7}')
$$
 and for any concave function $M$: $[0,\infty)\to {\mathbb R}$ such that $M(0)=0$,  the following inequality is true:
$$
 \sum_{k=1}^\infty p_k b_k M(a_k )\ge \inf\limits_{s\in [1,\infty)}\left\{   M\bigg(\frac {\sum_{k=1}^\infty p_k a_k }{\sum_{k=1}^s p_k}\bigg) \sum_{k=1}^s p_k b_k\right\}.
 \eqno(\ref{206}')
$$}


\small

\begin{itemize}
  \item[{[1]}]
 D.\,S. Mitrinovi\'c, J.\,E. Pe\v cari\'c, A.\,M. Fink, "Classical and new inequalities in analysis."\ Kluwer (1993),
 740~p.
  \item[{[2]}]
 P.\,L. Chebyshev, "O priblizhennyh vyrazhenijah odnih integralov cherez drugie,"\ Soobschenija i Protokoly Zasedanij Matematicheskogo Obschestva pri Imperatorskom Khar'kovskom Universitete, No.~2 (1882), \mbox{93--98}; Polnoe Sobranie Sochinenii P.L.~Chebysheva. Moskva--Leningrad, 1978. {\bf ~3}, 128--131.
  \item[{[3]}]
 P.\,L. Chebyshev, "Ob odnom rjade, dostavljajuschem predel'nye velichiny integralov pri razlozhenii podintegral'noi funkcii na mnozhiteli,"\ Prilozhenie k 57 tomy Zapisok Imp. Akad. Nauk, No~4 (1883);  Polnoe Sobranie Sochinenii P.L.~Chebysheva. Moskva--Leningrad, 1978. {\bf ~3}, 157--169.
   \item[{[4]}]
D.\,S. Mitrinovi\'c and P.\,M. Vasi\'c,\
 "History, variations and generalisations of the \v Ceby\v sev inequality and the question of some priorities,"\ Univ. Beograd. Publ. Elektrotehn. Fak. Ser. Mat. Fiz. No.~{\bf 461--497} (1974),~1--30.
  \item[{[5]}]
D.\,S. Mitrinovi\'c and J.\,E. Pe\v cari\'c, "History, variations and generalisations of the \v Ceby\v sev inequality and the question of some priorities II,"\  Rad JAZU (Zagreb), {\bf 450}, fasc. 9 (1990), 139--156.
  \item[{[6]}]
H.\,P. Heinig and L. Maligranda, "Chebyshev inequality in function spaces", Real Anal. Exchange, {\bf 17},   1 (1991/92), 211–247.
 \item[{[7]}]
A.L. Shidlich, "On necessary and sufficient for validity of some Chebyshev-Type inequalities,"\ Journal of
Mathematical Inequalities,  {\bf 5},  1 (2011), 71–85.
  \item[{[8]}]
J.F. Steffensen, "En Ulighted mellem Middelve edier,"\  Mat. Tidsskrift, {\bf 1920 B}, 49--53.

\end{itemize}
\enddocument